\magnification=\magstep1 
\baselineskip=14pt 
\def\mapright#1{\smash{\mathop{\longrightarrow}\limits^{#1}}} 

% Reference macros. Sets up a count as a counter for refs, resets it to 0.
\newcount\refno\refno=0
% References are typed by putting \ref at the start of the line, the text
% of the reference, and a blank line afterwards.  After the blank line
% following a reference to (say) Ladner, put
%     \newcount\ladner\ladner=\refno
% Then in your document, type [\the\ladner] to obtain [XX]
% At the end of the list of references, the line :
%            \immediate\closeout\reffile
% should appear to close the references file, and then
%            \input TempReferences
% will read them in later.
\chardef\other=12
\newwrite\reffile
\immediate\openout\reffile=TempReferences
\outer\def\ref{\par\medbreak\global\advance\refno by 1
  \immediate\write\reffile{}
  \immediate\write\reffile{\noexpand\item{[\the\refno]}}
  \copytoblankline}
\def\copytoblankline{\begingroup\setupcopy\copyref}
\def\setupcopy{\def\do##1{\catcode`##1=\other}\dospecials
  \catcode`\|=\other \obeylines}
{\obeylines \gdef\copyref#1
  {\def\next{#1}%
  \ifx\next\empty\let\next=\endgroup %
  \else\immediate\write\reffile{\next} \let\next=\copyref\fi\next}}

%\ref ALLEN, S.D.; SINCLAIR, A.M.; SMITH, R.R., The ideal structure of the Haagerup tensor product of C$^*$-algebras, {\sl J. reine angew. Math.}, 442 (1993) 111-148.

%\newcount\ASS\ASS=\refno

%\ref ARCHBOLD, R.J., Topologies for primal ideals, {\sl J. London Math. Soc.}, (2) 36 (1987) 524-542.

%\newcount\RJA\RJA=\refno

%\ref ARCHBOLD, R.J.; BATTY, C.J.K., On factorial states of operator algebras, III, {\sl J. Operator Theory}, 15 (1986) 53-81.

%\newcount\AB\AB=\refno

%\ref ARCHBOLD, R.J.; KANIUTH, E.; SCHLICHTING, G.; SOMERSET, D.W.B., Ideal spaces of the Haagerup tensor product of C$^*$-algebras, {\sl Internat. J. Math.}, 8 (1997) 1-29.

%\newcount\AKSS\AKSS=\refno

%\ref ARCHBOLD, R.J.; SOMERSET, D.W.B., Quasi-standard C$^*$-algebras, {\sl Math. Proc. Camb. Phil. Soc.}, 107 (1990) 349-360.

%\newcount\AS\AS=\refno

\ref BECKHOFF, F., Topologies on the space of ideals of a Banach algebra, {\sl Stud. Math.}, 115 (1995) 189-205.

\newcount\Be\Be=\refno

\ref BECKHOFF, F., Topologies of compact families on the ideal space of a Banach algebra, {\sl Stud. Math.}, 118 (1996) 63-75.

\newcount\Bec\Bec=\refno

\ref BECKHOFF, F.,  Topologies on the ideal space of a Banach algebra and spectral synthesis, {\sl Proc. Amer. Math. Soc.}, 125 (1997) 2859-2866.

\newcount\Beck\Beck=\refno

%\ref BECKHOFF, F., An example, unpublished.

%\newcount\Beckh\Beckh=\refno

%\ref BROWDER, A., {\sl Introduction to Function Algebras}, W.A. Benjamin, New York, 1969.

%\newcount\Brow\Brow=\refno

%\ref COHN, P.M., {\sl Algebra, Vol. 3} (2nd edn.), Wiley, Chichester, 1991.

%\newcount\Coh\Coh=\refno

\ref COLE, B., {\sl One point parts and the peak point conjecture}, Ph.D. dissertation, Yale Univ., 1968.

\newcount\Cole\Cole=\refno

%\ref DALES, H.G.; McCLURE, J.P., Nonstandard ideals in radical convolution algebras on a half-line, {\sl Can. J. Math.}, 39 (1987) 309-321.

%\newcount\DM\DM=\refno

\ref DAVIE, A.M., Quotients of uniform algebras, {\sl J. London Math. Soc.}, (2) 7 (1973) 31-40.

\newcount\Dav\Dav=\refno

\ref DIEUDONN\'E, J., {\sl Foundations of Modern Analysis}, Academic Press, London, 1969.

\newcount\Dieu\Dieu=\refno

%\ref DIXON, P.G., Topologically irreducible representations and radicals in Banach algebras, {\sl Proc. London Math. Soc.}, (3) 74 (1997) 174-200.

%\newcount\Dixon\Dixon=\refno

%\ref DOMAR, Y., Extensions of the Titchmarsh convolution theorem with applications in the theory of invariant subspaces, {\sl Proc. London Math. Soc.}, 46 (1983) 288-300.

%\newcount\Dom\Dom=\refno

\ref FEINSTEIN, J.F., A non-trivial, strongly regular uniform algebra, {\sl J. London Math. Soc.}, (2) 45 (1992) 288-300.

\newcount\Fein\Fein=\refno

\ref FEINSTEIN, J.F., Strong regularity and characters for Banach function algebras, preprint, University of Nottingham, 1997.

\newcount\Feins\Feins=\refno

\ref FELL, J.M.G., A Hausdorff topology for the closed subsets of a locally compact non-Hausdorff space, {\sl Proc. Amer. Math. Soc.}, 13 (1962) 472-476.

\newcount\Fel\Fel=\refno

%\ref FELL, J.M.G., The dual spaces of Banach algebras, {\sl Trans. Amer. Math. Soc.}, 114 (1965) 227-250.

%\newcount\Fell\Fell=\refno

\ref GAMELIN, T.W., {\sl Uniform algebras}, Prentice-Hall, New Jersey, 1969.

\newcount\Gam\Gam=\refno

\ref KELLEY, J.L., {\sl General Topology}, Van Nostrand, Princeton, 1955.

\newcount\Kel\Kel=\refno

\ref KHANIN, L.G., Spectral synthesis of ideals in classical algebras of smooth functions, {\sl Functions Spaces, (Edwardsville, IL, 1994)}, 167-182, Lecture Notes in Pure and Appl. Math., 172, Dekker, New York, 1995.

\newcount\Kha\Kha=\refno

%\ref KITCHEN, J.W.; ROBBINS, D.A., Sectional representation of Banach modules, {\sl Pacific J. Math.}, 109 (1983) 135-156.

%\newcount\KR\KR=\refno

%\ref KITCHEN, J.W.; ROBBINS, D.A., Bundles of Banach algebras, {\sl Internat. J. Math. Math. Sci.}, 17 (1994) 671-680.

%\newcount\KRo\KRo=\refno

%\ref KITCHEN, J.W.; ROBBINS, D.A., Bundles of Banach algebras, II, {\sl Houston J. Math.}, 20 (1994) 435-451.

%\newcount\KRob\KRob=\refno

%\ref KOPPERMAN, R., Asymmetry and duality in topology, {\sl Gen. Top. and Appl.}, 66 (1995) 1-39.

%\newcount\Ko\Ko=\refno

%\ref MEYER, M.J., Continuous dense embeddings of strong Moore algebras, {\sl Proc. Amer. Math. Soc.}, 116 (1992) 727-735.

%\newcount\Mey\Mey=\refno

\ref MORTINI, R., Closed and prime ideals in the algebra of bounded analytic 
functions, {\sl Bull. Austral. Math. Soc.},  35 (1987) 213-229. 

\newcount\Mor\Mor=\refno

%\ref NEUMANN, M. M., Commutative Banach algebras and decomposable operators, {\sl Monatshefte Math.}, 113 (1992) 227-243.

%\newcount\Neu\Neu=\refno

\ref O'FARRELL, A.G., A regular uniform algebra with a continuous point 
derivation of infinite order, {\sl Bull. London Math. Soc.}, 11 (1979) 41-44.

\newcount\OFAR\OFAR=\refno

\ref SOMERSET, D.W.B., Ideal spaces of Banach algebras, {\sl Proc. London Math. Soc.}, to appear.

\newcount\Id\Id=\refno

%\ref SOMERSET, D.W.B., Spectral synthesis for Banach algebras, submitted.

%\newcount\Syn\Syn=\refno

\ref STOUT, E.L., {\sl The Theory of Uniform Algebras}, Bogden and Quigley, New York, 1971.

\newcount\Stout\Stout=\refno

%\ref THOMAS, M.P., A non-standard ideal of a radical Banach algebra of power series, {\sl Acta Math.}, 152 (1984) 199-217.

%\newcount\Tho\Tho=\refno

%\ref VASILIEV, N.B., C$^*$-algebras with finite dimensional irreducible representations, {\sl Uspekhi Mat. Nauk} 21 (1966) 135-154 (in Russian; English translation, {\sl Russian Math. Surv.}, 21, 137-155.)

%\newcount\Vas\Vas=\refno

\ref WILKEN, D.R., A note on strongly regular uniform algebras, {\sl Can. J. Math.}, 21 (1969) 912-914.

\newcount\Wilk\Wilk=\refno

%\ref WILLCOX, A.B., Some structure theorems for a class of Banach algebras, {\sl Pacific J. Math.}, 6 (1956) 177-192

%\newcount\Will\Will=\refno

\immediate\closeout\reffile

\centerline {\bf A NOTE ON IDEAL SPACES OF BANACH ALGEBRAS} 
\bigskip 
\centerline {\bf J. F. Feinstein and D. W. B. Somerset} 
\bigskip 
\footnote{ }{1991 Mathematics Subject Classification 46H10}
\bigskip 
\noindent {\bf Abstract} In a previous paper the second author introduced a 
compact topology $\tau_r$ on the space of closed ideals of a unital
Banach algebra 
$A$. If $A$ is separable then $\tau_r$ is either metrizable or else neither 
Hausdorff nor first countable. Here it is shown that $\tau_r$ is Hausdorff if 
$A$ is $C^1[0,1]$, 
but that if $A$ is a uniform algebra then $\tau_r$ is Hausdorff if and only if 
$A$ has spectral synthesis. An example is given of a strongly regular, uniform 
algebra for which every maximal ideal has a bounded approximate identity, but 
which does not have spectral synthesis. 
\bigskip 
\bigskip 
\noindent {\bf Introduction} 
\bigskip 
\noindent Following the work of Beckhoff [\the\Be], [\the\Bec], [\the\Beck], 
the second author introduced in [\the\Id] the topology $\tau_r$ on the lattice 
$Id(A)$ of closed, two-sided ideals of a unital Banach algebra $A$
(throughout this note, all Banach algebras will be assumed to be unital). 
This topology is 
defined using the quotient norms, and is invariant under a change to an 
equivalent norm on $A$. If $A$ is a commutative Banach algebra, then 
$\tau_r$ coincides, on the maximal ideal space of $A$, with the Gelfand 
topology, while if $A$ is a C$^*$-algebra $\tau_r$ coincides with the compact, 
Hausdorff topology $\tau_s$ on $Id(A)$. It was shown in [\the\Id] that 
$\tau_r$ 
is always compact, and that for separable Banach algebras a dichotomy occurs: 
either $\tau_r$ is metrizable, or else $\tau_r$ is neither Hausdorff nor first 
countable. The first possibility is known to occur for finite-dimensional 
Banach algebras, TAF-algebras, and Banach algebras with spectral synthesis, 
but 
the question of whether the second possibility ever occurs was left open. 

Two cases of particular interest are the disc algebra, and the algebra 
$C^1[0,1]$. In both cases a complete description of the set of closed
ideals is available, but 
there is some difficulty in describing the quotient norms. The purpose of this 
paper is to show that $\tau_r$ is Hausdorff for $C^1[0,1]$, 
but that if $A$ is a uniform algebra then $\tau_r$ is Hausdorff 
on $Id(A)$ if and only if 
$A$ has spectral synthesis. Thus for the disc algebra $\tau_r$ is neither 
Hausdorff nor first countable. An example is given of a strongly regular, 
uniform algebra for which every maximal ideal has a bounded approximate 
identity, but which does not have spectral synthesis. 

Let us now define the topology $\tau_r$, which is the join of two weaker 
topologies. The first is easily defined: $\tau_u$ is the weakest topology on 
$Id(A)$ for which all the norm functions $I\mapsto\Vert a+I\Vert$ $(a\in A,\ 
I\in Id(A))$ are upper semi-continuous. The other topology $\tau_n$ can be 
described in various different ways, but none is particularly easy to work 
with. A net $( I_{\alpha})$ in $Id(A)$ is said to have the {\sl normality 
property} with respect to $I\in Id(A)$ if $a\notin I$ implies that 
$\lim\inf\Vert a+I_{\alpha}\Vert >0$. Let $\tau_n$ be the topology whose 
closed 
sets $N$ have the property that if $(I_{\alpha})$ is a net in $N$ with the 
normality property relative to $I\in Id(A)$ then $I\in N$. It follows that if 
$(I_{\alpha})$ is a net in $Id(A)$ having the normality property relative to 
$I\in Id(A)$ then $I_{\alpha}\to I$ $(\tau_n)$. Any topology for which 
convergent nets have the normality property with respect to each of their 
limits (such a topology is said to have the {\sl normality property}) is 
necessarily stronger 
than $\tau_n$, but it is not clear that $\tau_n$ need have the normality 
property. Indeed the following is true. Let $\tau_r$ be the topology on 
$Id(A)$ 
generated by $\tau_u$ and $\tau_n$. Then $\tau_r$ is always compact [\the\Id; 
2.3], and $\tau_r$ is Hausdorff if and only if $\tau_n$ has the normality 
property [\the\Id; 2.12]. 
Since, as we shall show, $\tau_r$ is not Hausdorff for uniform algebras 
without 
spectral synthesis, it follows that $\tau_n$ does not have the normality 
property in these algebras. 

We shall need the following lemma. 
\bigskip 
\noindent {\bf Lemma 0.1} {\sl Let $(I_{\alpha})$ be a net in $Id(A)$, either 
decreasing or increasing, and correspondingly either set $I=\bigcap 
I_{\alpha}$ 
or $I=\overline{\bigcup I_{\alpha}}$. Then $I_{\alpha}\to I$ $(\tau_r)$.} 
\bigskip 
\noindent {\bf Proof.} In either case, it is easy to see that $(I_{\alpha})$ 
has the normality property with respect to $I$, and converges to $I$ 
$(\tau_u)$. Hence $I_{\alpha}\to I$ $(\tau_r)$. Q.E.D. 
\bigskip 
\bigskip 
\noindent {\bf 1. Uniform algebras} 
\bigskip 
\noindent In this section we show that if $A$ is a uniform algebra then 
$\tau_r$ is Hausdorff if and only if $A$ has spectral synthesis. An example is 
given of a strongly regular, uniform algebra in which every maximal ideal has 
a bounded approximate identity, but which does not have spectral synthesis. 

Let $A$ be a Banach function algebra, with Shilov boundary $\Gamma(A)$. 
For a non-empty, closed subset $F$ of $Max(A)$ let $I(F)$ be the ideal of 
elements of $A$ which vanish on $F$, and let $J(F)$ the ideal of elements 
vanishing in a neighbourhood of $F$ in $Max(A)$. If $F\subseteq\Gamma(A)$ let 
$L(F)$ be 
the ideal which is the closure of set of elements vanishing in a neighbourhood 
of $F$ in $\Gamma(A)$. It is not clear whether $L(F)$ need always be the 
closure of $J(F)$, but of course this is so whenever $Max(A)=\Gamma(A)$. 
Recall 
that $A$ has {\sl spectral synthesis} if $J(F)$ is dense in $I(F)$ for each 
non-empty, closed subset $F$ of $Max(A)$, and that $A$ is {\sl strongly 
regular} if $J(\{ x\})$ is dense in $I(\{ x\})$ for each $x\in Max(A)$. If 
$I(\{ x\})=L(\{ x\})$ for each $x\in \Gamma(A)$ then $A$ is {\sl strongly 
regular on $\Gamma(A)$}; in fact, this implies that $\Gamma(A)=Max(A)$, see 
[\the\Wilk; Lemma] (for uniform algebras) and [\the\Feins] (for the general 
case: the elementary argument used is due to Mortini, and may be found in 
[\the\Mor]). 

Now let $A$ be a uniform algebra. A subset $X$ of $Max(A)$ is a {\sl peak set} 
if there is an $f\in A$ 
such that $f(x)=1$ for all $x\in X$, and $\vert f(x)\vert<1$ for all $x\notin 
X$. An intersection of peak sets is called a {\sl 
p-set}, and a singleton p-set is a {\sl p-point}. We shall use the facts that 
the union of two p-sets is a p-set [\the\Gam; II.12.8], and the set of 
p-points 
is a dense subset of $\Gamma(A)$ [\the\Stout; 7.24]. 

\bigskip 
\noindent {\bf Lemma 1.1} {\sl Let $A$ be a uniform algebra, $X$ a closed 
subset of $Max(A)$, and $F$ a p-set in $Max(A)$. Suppose that $f\in I(X)$. 
Then 
$\Vert f+I(X\cup F)\Vert=\sup\{\vert f(y)\vert:y\in F\}$.} 
\bigskip 
\noindent {\bf Proof.} It is trivial that $\Vert f+I(X\cup 
F)\Vert\ge\sup\{\vert f(y)\vert:y\in F\}$. For the converse, there exists, by 
[\the\Stout; 7.22], an $h\in fA\subseteq 
I(X)$ such that $h|_F=f|_F$ and $\Vert h\Vert=\sup\{\vert f(y)\vert:y\in F\}$. 
Hence $h-f\in I(X\cup F)$, and the result follows. Q.E.D. 

\bigskip 
\noindent {\bf Theorem 1.2} {\sl Let $A$ be a uniform algebra. Then $\tau_r$ 
is 
Hausdorff on $Id(A)$ if and only if $A$ has spectral synthesis.} 
\bigskip 
\noindent {\bf Proof.} If $A$ has spectral synthesis then $\tau_r$ is 
Hausdorff 
on $Id(A)$ by [\the\Id; 3.1]. 

If $A$ does not have spectral synthesis, there are two possibilities. Either 
$Max(A)=\Gamma(A)$, in which case $J(F)$ is dense in $L(F)$ for every closed 
subset $F$ of $Max(A)$, so by assumption there is a closed subset $X$ of 
$Max(A)$ such that $I(X)\ne L(X)$. Otherwise $Max(A)\ne \Gamma (A)$, so by 
[\the\Wilk; Lemma] $A$ is not strongly regular on $\Gamma (A)$. In other 
words, 
there exists $x\in \Gamma (A)$ such that $I(\{ x\})\ne L(\{ x\})$. Thus in 
either event there is a closed subset $X$ of $\Gamma(A)$ such that $I(X)\ne 
L(X)$. 

Let $(N_{\alpha})_{\alpha}$ be a net of decreasing, closed neighbourhoods of 
$X$ in $\Gamma (A)$ such that $\bigcap N_{\alpha}=X$. Then $(I(N_{\alpha}))_{\alpha}$ is an 
increasing net in $Id(A)$, and $I(N_{\alpha})\subseteq L(X)$, for each 
$\alpha$, so $I:=\overline{\bigcup_{\alpha} I(N_{\alpha})}\subseteq L(X)$. 
For each $\alpha$, let $(F_{\beta(\alpha)})_{\beta(\alpha)}$ be the 
increasing net of p-sets in $N_{\alpha}$ (recall that the union of two p-sets 
is a p-set [\the\Gam; II.12.8]). Then 
$(I({F_{\beta(\alpha)}}))_{\beta(\alpha)}$ is a decreasing net in 
$Id(A)$, and the density of the set of p-points in $\Gamma (A)$ [\the\Stout; 
7.24] implies that $\bigcap_{\beta(\alpha)}I({F_{\beta(\alpha)}}) = 
I(N_{\alpha})$. Hence $I({F_{\beta(\alpha)}})\mapright{\beta(\alpha)} 
I(N_{\alpha})$ $(\tau_r)$ by Lemma 0.1. But $I(N_{\alpha})\to 
I$ $(\tau_r)$, also by Lemma 0.1, so if $(I(F_{\gamma}))_{\gamma}$ denotes the 
`diagonal' net, see [\the\Kel; \S 2, Theorem 4], then $I(F_{\gamma})\to I\subseteq L(X)$ $(\tau_r)$. 

Suppose that $f\notin I(X)$. Then there is an open subset $O$ of $\Gamma (A)$, 
meeting $X$, and an $\epsilon >0$ such that $\vert f(x)\vert>\epsilon$ for all 
$x\in O$. By the density of p-points in $\Gamma (A)$, there is, for each 
$\alpha$, a p-point in $O\cap N_{\alpha}$, so that there is a 
$\beta_0(\alpha)$ 
such that $\Vert f+I(F_{\beta(\alpha)})\Vert>\epsilon$ for all 
$\beta(\alpha)\ge\beta_0(\alpha)$. Hence the `diagonal' net 
$I(F_{\gamma})\to I(X)$ 
$(\tau_n)$. 
On the other hand, if $f\in I(X)$ then a simple topological argument shows 
that 
there exists $\alpha_0$ such that for $\alpha\ge\alpha_0$, 
$N_{\alpha}\subseteq\{x\in \Gamma(A):\vert f(x)\vert<\epsilon\}$. Thus for 
$\alpha\ge\alpha_0$, 
$$\Vert f+I(F_{\beta(\alpha)})\Vert\le \Vert f+I(X\cup 
F_{\beta(\alpha)})\Vert<\epsilon$$ 
for all $\beta(\alpha)$, by Lemma 1.1. Hence 
$I(F_{\gamma})\to I(X)$ $(\tau_u)$, 
using [\the\Id; 2.1], 
and so
$$I(F_{\gamma})\to I(X)~~(\tau_r).$$ 
Since $I\ne I(X)$, $\tau_r$ is not Hausdorff. Q.E.D. 
\bigskip 
\noindent In particular, $\tau_r$ is neither Hausdorff nor first countable if 
$A$ is the disc algebra. 

It seems to be an open question whether there are any uniform algebras with 
spectral synthesis, other than the trivial examples of commutative 
C$^*$-algebras. Examples are given in [\the\Fein] of non-trivial, strongly 
regular uniform algebras for which every maximal ideal has a bounded 
approximate identity. However, this is not enough to ensure spectral synthesis, 
as we shall now see. 
\smallskip 
We begin by restating some results from [\the\Fein], which use Cole's systems 
of root extensions [\the\Cole]. 
\bigskip 
\noindent {\bf Proposition 1.3} {\sl Let $A_1$ be a normal uniform algebra on 
a compact, Hausdorff space $X_1$. Then there is a uniform algebra $A$ on a 
compact, Hausdorff space 
$X$, a surjective continuous map $\pi$ from $X$ onto $X_1$ and a bounded 
linear map $S: A \rightarrow A_1$ with the following properties. 

(a) The uniform algebra $A$ is strongly regular and every maximal ideal 
in $A$ has a bounded approximate identity. 

(b) For every $f \in A_1$, $f \circ \pi$ is in $A$. 

(c) If $x\in X_1$ and $g \in A$ with $g$ constantly equal to some complex 
number $c$ on $\pi^{-1}(\{x\})$ then $(Sg)(x) = c$. In particular, for all 
$f\in A_1$, $S(f \circ \pi) = f$. 

\noindent If $X_1$ is metrizable, then in addition to the above properties we 
may also 
insist that $X$ is metrizable.} 
\bigskip 
\noindent 
Surprisingly, perhaps, these properties are not enough to guarantee that $A$ 
has spectral synthesis. Indeed, the following result shows that the method of 
[\the\Fein] 
cannot produce a uniform algebra with spectral synthesis unless it starts 
with such an algebra in the first place. 
\bigskip 
\noindent {\bf Proposition 1.4} {\sl Let $A_1$, $X_1$, $A$, and $X$ be as in 
Proposition 1.3. If $A$ has spectral synthesis, then $A_1$ has spectral 
synthesis.} 
\bigskip 
\noindent {\bf Proof.} Let $S$, $\pi$ be as in Proposition 1.3. 
Note that $A$ and $A_1$ are both normal, so that $X=Max(A)$, and 
$X_1=Max(A_1)$. 
Suppose that $A$ has spectral synthesis. 
Let $E$ be a closed subset of $X_1$. Let $f \in A_1$ with 
$f(E) \subseteq \{0\}$. We show that $f$ is a uniform limit of 
functions in $A_1$ each of which vanishes in a neighbourhood of $E$. 

Set $F=\pi^{-1}(E)$. Certainly $f\circ\pi$ vanishes on $F$. Since $A$ has 
spectral synthesis, there is a sequence of functions $g_n \in A$ each 
vanishing in a neighbourhood of $F$, and such that $g_n$ tends uniformly 
to $f\circ\pi$ on $X$. We then have that $S(g_n)$ tends uniformly to 
$f$ in $A_1$. Moreover it follows from Proposition 1.3(c), and some 
elementary topology, that each of the functions $S(g_n)$ vanishes in a 
neighbourhood of $E$. The result now follows. Q.E.D. 
\bigskip 
In fact the same argument shows, more generally, that it is impossible to 
obtain a uniform algebra with spectral synthesis from one without merely by 
taking a system of root extensions. 
\bigskip 
\noindent {\bf Theorem 1.5} {\sl There exists a strongly regular, uniform 
algebra $A$ on a compact, metric space such that every maximal ideal of 
$A$ has a bounded approximate identity, but such that $A$ does not have 
spectral synthesis.} 
\bigskip 
\noindent{\bf Proof.} 
Let $A_1$ be a normal, uniform algebra on a compact, metric space $X_1$ such 
that $A_1$ is not strongly regular (the uniform algebra constructed by 
O'Farrell in [\the\OFAR] would do). Certainly $A_1$ does not have spectral 
synthesis. The result now follows immediately from Propositions 1.3 and 1.4. 
Q.E.D. 
\bigskip 
\bigskip 
\noindent {\bf 2. The Banach algebra $C^1[0,1]$} 
\bigskip 
\noindent In this section we show that if $A$ is the Banach algebra $C^1[0,1]$ 
then $\tau_r$ is Hausdorff on $Id(A)$. It is interesting that it is possible 
to 
identify $\tau_r$ without having a complete description of the quotient norms 
of $C^1[0,1]$. 
\bigskip 
Let $A=C^1[0,1]$ be the Banach algebra of continuously differentiable, complex 
functions on $[0,1]$, with norm given by $\Vert f\Vert =\Vert 
f\Vert_{\infty}+\Vert f'\Vert_{\infty}$ where $\Vert\, .\,\Vert_{\infty}$ is 
the supremum norm, and $f'$ is the derivative of $f$. 
It is well known that each closed ideal $I$ of $A$ can be described uniquely 
in 
the form $$I=\{ f\in A: f(x)=0\ (x\in C),\ f'(x)=0\ (x\in E)\},$$ where $C$ is 
a closed subset of $[0,1]$ and $E$ is a closed subset of $C$ containing the 
limit points of $C$. For topological reasons, however, we shall use a slightly 
different description. 

For $(x,y)\in [0,1]\times [0,1]$ define a functional $\phi_{x,y}$ by 
$$\phi_{x,y}(f)={f(x)-f(y)\over x-y}$$ if $x\ne y$, and by 
$\phi_{x,y}(f)=f'(x)$ if $x=y$. With $C$ and $E$ as above, let $D$ be the 
closed subset of $[0,1]^2$ given by $$D=\{ (x,y)\in C\times C: {\rm either }\ 
x\ne y,\ {\rm or }\ x\in E\}.$$ It is straightforward to check that the ideal 
$I$ above can now be described as $$I=\{ f\in A: f(x)=0\ (x\in C),\ 
\phi_{x,y}(f)=0\ ((x,y)\in D)\}.$$ Furthermore, the pair $(C,D)$ of closed 
sets 
has the following two properties: (i) $D\subseteq C\times C$, (ii) $D$ 
contains 
all the non-diagonal points of $C\times C$. On the other hand, if $C$ and $D$ 
are any closed subsets of $[0,1]$ and $[0,1]^2$, respectively, having these 
two 
properties, then there is a unique ideal $I(C,D)\in Id(A)$ given by 
$$I(C,D)=\{ f\in A: f(x)=0\ (x\in C),\ \phi_{x,y}(f)=0\ ((x,y)\in D)\}.$$ 
For $Y=[0,1]$ or $Y=[0,1]^2$, let $H(Y)$ be the set of closed subsets of $Y$, 
equipped with the Fell topology (the topology induced by the Hausdorff metric) 
which is compact and Hausdorff [\the\Fel]. Let $K$ be the subset of 
$H([0,1])\times H([0,1]^2)$ consisting of pairs of closed sets $(C,D)$ having 
the properties (i) and (ii). By the remarks above, there is a bijective 
correspondence between $K$ and $Id(A)$. This is the description of $Id(A)$ 
which we shall use. 

\bigskip 
\noindent {\bf Lemma 2.1} {\sl $K$ is a closed subset of $H([0,1])\times 
H([0,1]^2)$, with the product Fell topology.} 
\bigskip 
\noindent {\bf Proof.} Let $((C_{\alpha},D_{\alpha}))_{\alpha}$ be a net in 
$K$ 
with limit $(C,D)\in H([0,1])\times H([0,1]^2)$. If $(x,y)\in D$, there is a 
net $((x_{\alpha}, y_{\alpha}))_{\alpha}$ with $(x_{\alpha},y_{\alpha})\in 
D_{\alpha}$ for each $\alpha$, such that $(x_{\alpha},y_{\alpha})\to (x,y)$. 
By 
property (i), $x_{\alpha},y_{\alpha}\in C_{\alpha}$ for each $\alpha$, and 
$x_{\alpha}\to x$, $y_{\alpha}\to y$, so $x,y\in C$. Hence $D\subseteq C\times 
C$. 

On the other hand, suppose that $x,y\in C$ with 
$x\ne y$. Then there are nets $(x_{\alpha})$ and $(y_{\alpha})$ with 
$x_{\alpha},y_{\alpha}\in C_{\alpha}$ for each $\alpha$ and $x_{\alpha}\to x$, 
$y_{\alpha}\to y$. By passing to subnets we may assume that $x_{\alpha}\ne 
y_{\alpha}$ for each $\alpha$. Hence $(x_{\alpha},y_{\alpha})\in D_{\alpha}$ 
by 
property (ii), and $(x_{\alpha},y_{\alpha})\to (x,y)$, so $(x,y)\in D$. 

This establishes that (i) and (ii) hold for the pair $(C,D)$, so $(C,D)\in K$. 
Q.E.D. 
\bigskip 
\noindent Thus $K$ is a compact, Hausdorff space in the restriction of the 
product Fell topology. Let $\tau$ be the corresponding compact, Hausdorff 
topology on $Id(A)$, induced by the bijective correspondence 
$(C,D)\leftrightarrow I(C,D)$ between $K$ and $Id(A)$. We shall show that 
$\tau$ coincides with $\tau_r$. 

Let $I=I(C,D)$ be a closed ideal of $A$, and define a semi-norm $r_I$ on $A$ 
by 
$$r_I(f)=\sup\{\vert f(x)\vert:x\in C\}+\sup\{\vert 
\phi_{x,y}(f)\vert:(x,y)\in 
D\}\ \ \ \ \ (f\in A).$$ Note that $\ker r_I=I$, and that if $I_{\alpha}\to I$ 
$(\tau)$ then $r_{I_{\alpha}}(f)\to r_I(f)$ for all $f\in A$. Let $q_I$ be the 
quotient semi-norm $q_I(f)=\Vert f+I\Vert$. 
\bigskip 
\noindent {\bf Lemma 2.2} {\sl Let $I=I(C,D)$ be a closed ideal of $A$, and 
let 
$f\in A$. Then $q_I(f)\ge r_I(f)$. } 
\bigskip 
\noindent {\bf Proof.} Let $h\in A$ with $h-f\in I$. Then $(h-f)(x)=0$, for 
$x\in C$, and $\phi_{x,y}(h-f)=0$, for $(x,y)\in D$. Hence $\Vert 
h\Vert_{\infty}\ge\sup\{ \vert 
f(x)\vert: x\in C\}$, by common sense, and $\Vert h'\Vert_{\infty}\ge 
\sup\{\vert \phi_{x,y}(f)\vert:(x,y)\in D\}$, by the Mean Value Theorem 
[\the\Dieu; 8.5.4] (recall that we are dealing with complex functions). Thus 
$\Vert h\Vert\ge r_I(f)$. Q.E.D. 
\bigskip 
\noindent Let $S_1(A)$ denote the set of semi-norms $\rho$ on $A$ satisfying 
$\rho(f)\le \Vert f\Vert$ for all $f\in A$. Then $S_1(A)$ is a compact, 
Hausdorff space with respect to the topology of pointwise convergence 
[\the\Be]. Note that if $(\rho_{\alpha})$ is a net in $S_1(A)$ and 
$\rho_\alpha\to\rho\in S_1(A)$ then $\ker\rho_{\alpha}\to\ker\rho$ $(\tau_u)$. 
\bigskip 
\noindent {\bf Theorem 2.3} {\sl The topologies $\tau$ and $\tau_r$ coincide 
on 
$Id(A)$. Hence $\tau_r$ is a compact, Hausdorff topology.} 
\bigskip 
\noindent {\bf Proof.} First we show that $\tau$ has the normality property. 
Suppose that $I_{\alpha}=I_{\alpha}(C_{\alpha}, D_{\alpha})\to I=I(C,D)$ 
$(\tau)$, and that $f\notin I$. If there exists $x\in C$ such that $\vert 
f(x)\vert=\delta>0$, then eventually the set $\{ y\in [0,1]:\vert f(y)\vert 
>\delta/2\}\cap C_{\alpha}$ is non-empty, so eventually $\Vert 
f+I_{\alpha}\Vert\ge r_{I_{\alpha}}(f)>\delta/2$ by Lemma 2.2. Otherwise there 
exists $(x,y)\in D$ such that $\vert\phi_{x,y}(f)\vert =\delta>0$, so 
eventually the set $\{ (z,w)\in [0,1]^2:\vert \phi_{z,w}(f)\vert 
>\delta/2\}\cap D_{\alpha}$ is non-empty, so eventually $\Vert 
f+I_{\alpha}\Vert>\delta/2$, again by Lemma 2.2. Thus we have shown that 
$\tau$ 
has the normality property for the net $(I_{\alpha})$ with respect to $I$. 

Next we show that $\tau\ge\tau_u$. Suppose, for a contradiction, that 
$I_{\alpha}=I_{\alpha} 
(C_{\alpha}, D_{\alpha})\to I=I(C,D)$ $(\tau)$, but that $I_{\alpha}\not\to I$ 
$(\tau_u)$. We may assume, by passing to a subnet, if necessary, and using the 
compactness of $S_1(A)$, that $(q_{I_{\alpha}})$ is a convergent net in 
$S_1(A)$, with limit $\sigma$, say, and that $I\ne\ker\sigma$. Set 
$J=I(F,G)=\ker\sigma$, and note that since $q_{I_{\alpha}}(f)\ge 
r_{I_{\alpha}}(f)$ ($f\in A$) for all $\alpha$, it follows that 
$\sigma(f)=\lim 
q_{I_{\alpha}}(f)\ge \lim r_{I_{\alpha}}(f)=r_I(f)$. Hence $J\subseteq I$, so 
$F\supseteq C$ and $G\supseteq D$. Suppose that $x\in F\setminus C$, and let 
$f\in A$ be such that $f(x)=1$, and $f$ vanishes in an open neighbourhood $N$ 
of $C$. Then $D\subseteq C\times C\subseteq N\times N$, so $\phi_{y,z}(f)=0$ 
for 
all $(y,z)\in D$. Hence $f\in I$, but $f\notin J$ so $\sigma(f)>0$. Since 
$I_{\alpha}\to I$ 
$(\tau)$, eventually $C_{\alpha}\subseteq N$, so eventually $f\in I_{\alpha}$. 
Hence $0=\lim q_{I_{\alpha}}(f)=\sigma(f)>0$, a contradiction. Thus $F=C$. 

Now suppose that $(x,y)\in G\setminus D$. Then $x=y$, by property (ii), since 
$C=F$, so $x$ is an isolated point of $C$, by the first description of 
$Id(A)$. 
Let $f\in A$ such that $f(x)=0$, $f'(x)\ne 0$, and $f$ vanishes on a closed 
set 
$M$, disjoint from $x$, but containing $C\setminus\{ x\}$ in its interior. 
Thus 
$f\in I\setminus J$. Since $I_{\alpha}\to I$ $(\tau)$, there exists $\alpha_0$ 
such that for $\alpha\ge\alpha_0$, $C_{\alpha}\cap ([0,1]\setminus M)$ is a 
singleton $\{ x_{\alpha}\}$. For fixed $\alpha\ge\alpha_0$ and any $\epsilon 
>0$ there exists $g\in A$ such that $g(x_{\alpha})=f(x_{\alpha})$, $g(M)=0$, 
and $$\Vert g\Vert<\vert f(x_{\alpha})\vert+ {\vert f(x_{\alpha})\vert\over 
d(x_{\alpha},M)}+\epsilon,$$ where $d(x_{\alpha},M)$ is the distance from 
$x_{\alpha}$ to the set $M$. Since $\epsilon$ is arbitrary and $g-f\in 
I_{\alpha}$ it follows that $$q_{I_{\alpha}}(f)\le \vert f(x_{\alpha})\vert+ 
{\vert f(x_{\alpha})\vert\over d(x_{\alpha},M)}.$$ But $f(x_{\alpha})\to 
f(x)=0$, and $d(x_{\alpha},M)\to d(x,M)>0$, so $\sigma(f)=\lim 
q_{I_{\alpha}}(f)=0$. Hence $f\in J$, a contradiction. It follows that $G=D$, 
so $I=J$, another contradiction. Thus $I_{\alpha}\to I$ $(\tau_u)$ after all. 

We have shown that $\tau$ is a compact (Hausdorff) topology stronger than 
$\tau_u$, and 
such that $\tau$-convergent nets have the normality property. Hence 
$\tau=\tau_r$ by [\the\Id; 2.11]. Q.E.D. 
\bigskip 
\noindent {\bf Remarks.} (i) $C^1[0,1]$ is, with an equivalent norm, 
isometrically isomorphic to a quotient of a uniform algebra [\the\Dav]. It was 
shown in [\the\Id; 2.9] that if $I$ is a closed ideal of a Banach algebra $A$ 
then the natural map from $Id(A/I)$ onto the set $\{ J\in Id(A): J\supseteq 
I\}$ is a $\tau_r$--$\tau_r$ homeomorphism. 

(ii) Khanin [\the\Kha] has introduced the class of D-algebras, of which 
$C^1[0,1]$ and commutative Banach algebras with spectral synthesis are 
examples. It would be interesting to know if $\tau_r$ is Hausdorff on $Id(A)$ 
for all D-algebras $A$. 
\bigskip 
\bigskip 
\centerline {\bf References} 
\medskip 
\input TempReferences.tex 
\bigskip 
\centerline{Department of Mathematics} 
\centerline{University of Nottingham} 
\centerline{NG7 2RD} 
\centerline{U.K.} 
\smallskip 
\centerline{e-mail: jff@maths.nott.ac.uk} 
\bigskip 
\centerline{Department of Mathematical Sciences} 
\centerline{University of Aberdeen} 
\centerline{AB24 3UE} 
\centerline{U.K.} 
\smallskip 
\centerline{e-mail: ds@maths.abdn.ac.uk} 

\end